\documentstyle[12 pt]{article}
\def \r  {{\bf{R}}}
\def \om {\Omega}
\def \ra {\rightarrow}

\def \p  {\partial}
\def \l  {\lambda}
\def\endpf{\hbox{\vrule height1.5ex width.5em}}
\def \d  {\displaystyle}
\def \f  {\frac}
\def \e  {\varepsilon}

\def\xe  {\chi_{E}}
\def\xw  {\chi_{\omega}}
\def\w   {\widetilde}
\def\uad {{\mathcal U}_{ad}}

\begin{document}
\title{A  Bang-Bang Principle of Time Optimal Internal Controls of the Heat
Equation  \footnote{ This work was supported by the National Natural
Science Foundation of China under Grants 60574071 and 10471053, and
by the key project of Chinese Ministry of Education. }}
\author{Gengsheng Wang \\
 School of Mathematics and Statistics, Wuhan University,\\ Wuhan, Hubei,
430072, P.R.of China\\
\verb "wanggs@public.wh.hb.cn"}
\date{}
\maketitle
{\bf Abstract.} In this paper, we study a time optimal internal
control problem governed by the heat equation in $\Omega\times
[0,\infty)$. In the problem, the target set $S$ is nonempty in
$L^2(\Omega)$,  the control set $U$  is
 closed, bounded and  nonempty in  $L^2(\Omega)$ and control functions
  are taken from the set
  $\uad=\{u(\cdot, t): [0,
\infty)\ra L^2(\Omega)\; \mbox{measurable}; u(\cdot, t)\in U,
\;\mbox{a.e. in t}\;\}$.  We first establish a certain
 null controllability for the heat equation in $\Omega\times
[0,T]$, with controls restricted to a product set of an open
nonempty subset in $\Omega$ and a subset  of positive measure in the
interval $[0,T]$. Based on this, we prove that  each optimal control
$u^*(\cdot, t)$ of the problem satisfies necessarily the bang-bang
property: $u^*(\cdot, t)\in \p U$ for almost all $t\in [0, T^*]$,
where $\p U$ denotes the boundary of the set $U$ and $T^*$ is the
optimal time. We also obtain the uniqueness of the optimal control
when the target set $S$ is convex and the control set $U$ is a
closed ball.


\vskip 1cm
\section{Introduction}

\hspace*{0.5 cm} Let $\om$ be a bounded domain in  $\r^n, n\geq 1$,
with a $C^\infty$-smooth boundary. Let $\omega$ be an open subset of
$\om$. Denote by  $\xw$ the characteristic function of $\omega$.
Consider the following controlled heat equation:
$$
\left\{\begin{array}{ll}
y_t(x,t) -\Delta
y(x,t)=\xw(x)u(x,t)\;\;\;\;\;\;\;\;\;\;\;\;\;\;\;&
\mbox{in }\;\;\om\times(0,\infty),\\
y(x,t)=0 & \mbox{on }\;\; \p\om\times(0,\infty),\\
y(x,0)=y_0(x)&\mbox{in }\;\;\om,
\end{array}\right.\eqno{(1.1)}
$$
where $y_0(\cdot)$ is a function in $L^2(\om)$ and $u(x,t)$ is a
control function taken from the set of functions as follows:
$$
{\mathcal{U}}_{ad}=\{v: [0,\infty)\ra L^2(\om)\;\mbox{measurable;
}v(\cdot, t)\in U\;\mbox{for almost all }t\geq 0 \}.\eqno{(1.2)}
$$
Here $U$ is a closed, bounded and nonempty subset in $L^2(\Omega)$.
Notice that the control function $u$ is acted internally (or
locally) into the equation (1.1). If $\omega=\om$, we say that the
control is acted globally into the equation. We shall denote by
$y(x,t;u,y_0)$ or $y(x,t)$ the solution of the equation (1.1) if
there is no risk of causing confusion.

\medskip
In this paper, we shall study the following time optimal control
problem:
$$
 \bf{(P)}\;\; \;\; \; \mbox{Inf }\{\w t;\;y(\cdot,\w t;u,y_0)\in S,\;u\in
{\mathcal U}_{ad}\}.
$$
Where  $S$ is a nonempty subset in $L^2(\Omega)$. We call the set
$S$ as the target set, the set $U$ as the control set, the set
$\uad$ as the  control function set and $y_0$ as the initial state
for the problem $\bf (P)$. For simplicity, we shall call a control
function as a control.
 The number
 $$
 T^*\equiv \mbox{Inf }\{\w t;\;y(x,\w
t; u,y_0)\in S,\;u\in{\mathcal U}_{ad}\}
$$
 is called the optimal time
for the problem ${\bf (P)}$, a control $u^*\in\uad$ having the
property:
$$
y(x,T^*; u^*,y_0)\in S,
$$
is called an optimal control (or a time optimal control) for the
problem $\bf (P)$, and a control $u\in \uad$ having the property:
$$
y(x,T; u,y_0)\in S\;\;\;\mbox{for a certain positive number}\; T,
$$
 is called an admissible control for the problem $\bf (P)$.

 In this paper, we obtain that each optimal control
$u^*$ for the problem $\bf (P)$ satisfies the  bang-bang property:
$u^*(\cdot,t) \in \p U$ for almost all $t\in[0,T^*]$. We further
show that if the control set $U$ is a closed ball $B(0,R)$, centered
at the origin of $L^2(\Omega)$ and of  positive radius $R$, then
each optimal control $u^*$ for the problem $\bf (P)$ satisfies the
property: $\|\chi_\omega u^*(\cdot,t)\|_{L^2(\Omega)}=R$ for almost
all $t\in [0, T^*]$.
 We also prove
the uniqueness of the optimal control for the problem $\bf (P)$,
when the target set $S$ is  convex and nonempty and the control set
$U$ is a closed ball. Combining these with the existence result of
time optimal controls obtained in [17], ( See also [14].) we derive
that if the target set $S$ is a closed, convex and nonempty subset,
which contains the origin of $L^2(\Omega)$, and if the control set
$U$ is the ball $B(0,R)$, then the problem $\bf (P)$ has a unique
optimal control $u^*$ satisfying the bang-bang property:
$\|\chi_\omega u^*(\cdot, t)\|_{L^2(\Omega)}=R$ for almost all $t\in
[0, T^*]$.

 The bang-bang
principle above can be explained physically as follows: If an
outside force $u^*$, acted in an open subset $\omega$ of $\om$ and
with the maximum norm bound: $\|u^*(\cdot,t)\|_{L^2(\om)}\leq R$ for
almost all $t$,  makes the temperature distribution in $\om$ change
from an initial distribution $y_0(x)$ into the target set $S$ in the
shortest time $T^*$, then  $u^*$  takes necessarily the maximum norm
for almost all $t$ in $[0,T^*]$, namely, $\|\chi_\omega u^*(\cdot,
t)\|_{L^2(\Omega)}=R$ for almost all $t\in [0, T^*]$. This bang-bang
principle  is a weaker form if it is compared with the following
stronger form: {\it If $u^*$ is an
 optimal control of the problem
$\bf (P)$ where the  control function  set is
$$
\{ u(x,t)\in L^\infty(\Omega\times [0, \infty));\; |u(x,t)|\leq
R\;\; \mbox{for almost all} \; (x,t)\},
$$
 then $|u^*(x,t)|=R$ for
almost all $(x,t)\in\Omega\times [0,T^*]$.}

In this work, we observe that the bang-bang principle for the
problem $\bf (P)$ is based on the following  null controllability
property for the heat equation:

\medskip
${\bf (C) } $ {\it Let $T$ be a positive number and  let $E$ be a
subset of positive measure in the interval $[0,T]$.  For each
$\delta\geq 0$, we write  $E_\delta$ for the set $\{ t\in \r^1;
t+\delta\in E\}$ and denote by $\chi_{E_\delta}$ the characteristic
function of the set $E_\delta$. Then there exists a number
$\delta_0$ with $0< \delta_0 <T$ such that for each $\delta$ with
$0\leq \delta\leq \delta_0$ and for each element $y_0$ in
$L^2(\Omega)$, there is a control $u_\delta$ in the space
$L^\infty(0,T-\delta;L^2(\om))$ such that the solution $z^\delta$ to
the following controlled heat equation:
$$
\left\{\begin{array}{ll} z^\delta_t(x,t) -\Delta
z^\delta(x,t)=\chi_{E_\delta}(t)\xw(x)u_\delta(x,t)\;\;\;\;\;\;\;\;\;\;\;\;\;\;\;&
\mbox{in }\;\;\om\times(0,T-\delta),\\
z^\delta(x,t)=0 & \mbox{on }\;\; \p\om\times(0,T-\delta),\\
z^\delta(x,0)=y_0(x)&\mbox{in }\;\;\om,
\end{array}\right.\eqno{}
$$
satisfies $z^\delta(x,T-\delta)=0$ over $\Omega$. Moreover, the
control $u_\delta$ satisfies the following estimate:
$$
\|u_\delta\|^2_{L^\infty(0,T-\delta;L^2(\Omega))}\leq
L\|y_0\|^2_{L^2(\Omega)},
$$
where $L$ is a positive number independent of $\delta$ and $y_0$.}
\medskip

It is well known  that the null controllability $\bf (C)$ is
equivalent to the following observability inequality:

\medskip
$\bf (O) $ {\it There exist  positive numbers  $L$ and $\delta_0$
with $\delta_0<T$ such that
$$
[\d\int_{\om}(p^\delta(x,0))^2dx]^{\f{1}{2}}\leq
L\int_0^{T-\delta}\{\int_{\om}[\chi_{E_\delta}(t)\xw(x)p^\delta(x,t)]^2dx\}^{\f{1}{2}}dt
$$
for each number $\delta$ with $0\leq \delta \leq \delta_0$ and each
function $p^\delta_T(x)\in L^2(\om)$. Where $p^\delta(x,t)$ is the
solution to the following adjoint equation:
$$
\left\{\begin{array}{ll} p^\delta_t(x,t) +\Delta
p^\delta(x,t)=0\;\;\;\;\;\;\;\;\;\;\;\;\;\;\;&
\mbox{in }\;\;\om\times(0,T-\delta),\\
p^\delta(x,t)=0 & \mbox{on }\;\; \p\om\times(0,T-\delta),\\
p^\delta(x,T-\delta)=p^\delta_T(x)&\mbox{in }\;\;\om.
\end{array}\right.\eqno{}
$$}However the inequality $\bf (O)$ is not a trivial consequence of the
Carleman inequality for linear parabolic equation given in [6]. We
establish the property $\bf (C)$  by applying an iterative argument
stimulated by that in [7]. (See also [8] and [12].) Our iterative
argument is based on a sharp observability estimate on the
eigenfunctions of the Laplacian, due to G. Lebeau and E. Zuazua in
[8] ( See also [7].) and a special result in the measure theory
given in [13].
\medskip

It should be mentioned  that the problem $\bf (P)$ may have no
admissible control in many cases. For instance, if the target set
$S$ is a closed ball $ B(y_1, R)$  in $L^2(\Omega)$,  centered at
$y_1$ and of positive radius $R$ and if the control set $U$ is the
closed ball $B(0,1)$ in $L^2(\Omega)$, centered at the origin  and
of radius $1$, then a necessary condition for the existence of an
admissible control for the  problem $\bf (P)$ is as follows: ( See
[14].)
 $$
\|y_1\|_{L^2(\Omega)}\leq (1+\d\f{1}{{\l}_1}) (\|y_0\|_{L^2(\Omega)}
+1) +R,
$$
 where ${\l}_1$ is the first eigenvalue of the Laplacian.
 However, it was proved in [17] ( See also [14].) that when the target set $S$ is the origin of
 $L^2(\Omega)$ and the control set
 is the ball $B(0,R)$ with $R>0$,
 then the problem $\bf (P)$ has at least one time optimal
 control. From this, it follows that if the target set $S$ is a closed
 and convex subset,  which contains the origin of $L^2(\Omega)$, and if the control set
 $U$ is the ball $B(0,R)$ with $R>0$, then the problem
 $\bf (P)$    has at least one optimal control.

\medskip
The time optimal control problems for parabolic equations have been
extensively studied in the past  years. Here, we mention the works
[14], [17], [18] and [19], where the existence of  time optimal
controls for linear and some semi-linear parabolic equations was
investigated. We mention the works  [10] and [20], where both the
existence and the maximum principle of  time optimal controls
governed by certain parabolic equations were studied.  We mention
the works [1] and [9], where the maximum principle for  time optimal
controls was derived. We mention the works [3], [4], [5] and [11],
where the bang-bang principle ( in the weaker form) for time optimal
controls governed by linear parabolic and hyperbolic equations with
the controls acted in the whole domain $\Omega$ or the whole
boundary $\p\Omega$  was established. We mention the work [16],
where the bang-bang principle (in the stronger form) of time optimal
controls for the heat equation where the control is restricted in
the whole boundary was obtained. We mention the work [13], where the
bang-bang principle (in the stronger form) for  time optimal
controls of the one-dimensional heat equation where the control is
restricted in one ending point of the one-dimensional state space,
was derived. Moreover, the authors in [13] observed  that such a
bang-bang principle is based on a certain exactly boundary
null-controllability for  the one-dimensional heat equation from
arbitrary sets of positive measure in the time variable space. We
also mention a more recent work [18], where the bang-bang principle
(in the weaker form) of time optimal internal controls governed by
the heat equation  and with a ball centered at $0\in L^2(\Omega)$
and of a positive radius as the target was obtained. Moreover, in
[18], the bang-bang principle was obtained by  a certain unique
continuation property for the heat equation involving a measurable
set, and the maximum principle for the optimal controls.

This paper is organized as follows. In Section 2, we establish the
null controllability {\bf (C)}. In Section 3, we give and prove the
main results of the paper, namely, the bang-bang principle and the
uniqueness of the optimal control  for the problem $ \bf (P)$.
\section{The  null controllability $\bf (C)$}

\hspace*{0.5 cm} Let $T$ be a positive number and $E$ be a subset of
positive measure in the interval $[0,T]$. We denote by $m(E)$  the
Lebesgue measure of the set $E$ in $\r^1$. For each $\delta\geq 0$,
we write $E_\delta$ for the set $\{ t\in \r^1 ; \;\;t+\delta \in
E\}$ and denote by $\chi_{E_\delta}$ the characteristic function of
the set $E_\delta$. In what follows, we shall omit $(x,t)$ (or $t$)
in functions of $(x,t)$ (or functions of $t$), if there is no risk
of causing confusion. For each positive number $\delta$, we consider
the following controlled equation:
$$
\left\{\begin{array}{ll} y_t(x,t) -\Delta
y(x,t)=\chi_{E_\delta}(t)\xw(x)u(x,t)\;\;\;\;\;\;\;\;\;\;\;\;\;\;\;&
\mbox{in }\;\;\om\times(0,T-\delta),\\
y(x,t)=0 & \mbox{on }\;\; \p\om\times(0,T-\delta),\\
y(x,0)=y_0(x)&\mbox{in }\;\;\om,
\end{array}\right.\eqno{(2.1)}
$$
where $y_0\in L^2(\om)$ is a given function. The main result of this
section is as follows:

\bigskip
{\bf Theorem 2.1}. {\it Let $T$ be a positive number and let  $E$ be
a subset of positive measure in the interval $[0,T]$. Then there
exists a positive number $\delta_0$ with $\delta_0<T$ such that  for
each number $\delta$ with $0\leq\delta\leq\delta_0$ and for each
element $y_0$ in the space $L^2(\om)$, there is a control $u_\delta$
in the space $L^\infty (0,T-\delta;L^2(\om))$ with the estimate
$$
\|u_\delta\|^2_{ L^\infty (0,T-\delta;L^2(\om))}\leq
L\|y_0\|^2_{L^2(\Omega)}
$$
 for  a certain positive
constant $L$ independent of $\delta$ and $y_0$,  such that the
solution $y^\delta(x,t)$ to the equation (2.1) with $u$ being
replaced by $u_\delta$ reaches zero value at time $T-\delta$,
namely, $ y^\delta(x,T-\delta)=0\;\mbox{over }\om. $}

\bigskip
The proof of Theorem 2.1 is based on a  sharp estimate on the
eigenfunctions of the Laplacian due to G. Lebeau and E. Zuazua ( See
[8].) and a fundamental result in the measure theory, which will be
given in the later. Let $\{\l_i\}^\infty_{i=1},\;0<\l_1 < \l_2 \leq
\cdots$, be the eigenvalues of $-\Delta$ with the Dirichlet boundary
condition and $\{X_i(x)\}^\infty_{i=1}$ be the corresponding
eigenfunctions, which serve as an orthonormal basis of $L^2(\om)$.
Then we have the following result. (See [8].)

\bigskip
{\bf Theorem 2.2.}  {\it There exist two positive constants
$C_1,\;C_2>0$ such that
$$
\sum_{\l_i\leq r}|a_i|^2\leq C_1
e^{C_2\sqrt{r}}\d\int_{\omega}|\sum_{\l_i\leq r}a_i X_i(x)|^2dx
$$
for every finite $r>0$ and every choice of the coefficients
$\{a_i\}_{\l_i\leq r}$ with $a_i\in \r^1$.}

\medskip
Now, we shall first use Theorem 2.2 to derive a certain
controllability result, which will help us in the proof of Theorem
2.1. For each $r>0$, we set $ {\bf X_r}=\mbox{span
}\{X_i(x)\}_{\l_i\leq r}$, and consider the following dual equation:
$$
\left\{\begin{array}{ll} \varphi_t(x,t) +\Delta
\varphi(x,t)=0\;\;\;\;\;\;\;\;\;\;\;\;\;\;\;&
\mbox{in }\;\;\om\times(0,T),\\
\varphi(x,t)=0 & \mbox{on }\;\; \p\om\times(0,T),\\
\varphi(x,T)\in {\bf X_r}.
\end{array}\right.\eqno{(2.2)}
$$
Here, each element  $\varphi (x,T)$ in ${\bf X_r}$ can be written as
$$
\varphi
(x,T)=\sum_{\l_i\leq r}a_i X_i(x),
$$
for a certain sequence of real numbers  $\{a_i\}_{\l_i\leq r}$. Then
the solution $\varphi(x,t)$ to the equation (2.2) can be expressed
by
$$
\varphi (x,t)=\sum_{\l_i\leq r}a_i e^{-\l_i
(T-t)}X_i(x)\;\;\mbox{for all }t\in [0,T].
$$
Set $b_i(t)=a_i e^{-\l_i (T-t)},\;t\in [0,T].$ Then by Theorem 2.2,
we have
$$
\begin{array}{ll}
\d\sum_{\l_i\leq r}|b_i(t)|^2&\leq C_1 e^{C_2 \sqrt r}\d\int_\omega
|\sum_{\l_i\leq r}b_i(t)X_i(x)|^2dx \\
&=C_1 e^{C_2 \sqrt r}\d\int_\omega |\varphi(x,t)|^2dx\;\;\mbox{for
all }t\in [0,T].
\end{array}
$$
On the other hand,
$$
\begin{array}{ll}
\d\sum_{\l_i\leq r}|b_i(t)|^2=\sum_{\l_i\leq r}a^2_i
e^{-2\l_i(T-t)}&\geq \d\sum_{\l_i\leq r}a^2_i e^{-2\l_i T}\\
&=\d\int_\om \varphi^2(x,0)dx\;\;\;\;\mbox{for all }t\in [0,T].
\end{array}
$$
Hence,
$$
\d\int_\om \varphi^2(x,0)dx\leq C_1 e^{C_2\sqrt r} \int_\omega
|\varphi(x,t)|^2dx\;\;\;\;\mbox{for all }t\in [0,T],
$$
or equivalently,
$$
[\d\int_\om \varphi^2(x,0)dx]^{\f{1}{2}}\leq (C_1 e^{C_2\sqrt
r})^{\f{1}{2}}[\int_\omega
|\varphi(x,t)|^2dx]^{\f{1}{2}}\;\;\mbox{for all}\; t\in [0,T],
$$
from which, it follows that
$$
\d\int_E[\d\int_\om \varphi^2(x,0)dx]^{\f{1}{2}}dt\leq (C_1
e^{C_2\sqrt r})^{\f{1}{2}}\int_E [\int_\omega
|\varphi(x,t)|^2dx]^{\f{1}{2}}dt.
$$
Namely, we obtained that for each $\varphi(\cdot,T)\in {\bf X_r}$,
$$
\begin{array}{ll}
\d\int_\om \varphi^2(x,0)dx&\leq \d\f{C_1 e^{C_2\sqrt
r}}{(m(E))^2}\{\int_0^T[\int_\Omega
|\xe(t)\xw(x)\varphi(x,t)|^2dx]^{\f{1}{2}}dt\}^{2}\\
&=\d\f{C_1 e^{C_2\sqrt
r}}{(m(E))^2}\|\xe\xw\varphi\|^2_{L^1(0,T;L^2(\om))}.
\end{array}
\eqno{(2.3)}
$$
Write  $P_r$ for the orthogonal projection from $L^2(\om)$ to ${\bf
X_r}$. We next use (2.3) to obtain the following controllability
result.

\bigskip
{\bf Lemma 2.3}. {\it For each $r>0$, there exists a control $u_r$
in the space $L^\infty(0,T;L^2(\om))$ with the estimate
$$
\|u_r\|_{L^\infty(0,T;L^2(\om))}\leq \d\f{C_1 e^{C_2\sqrt
r}}{(m(E))^2}\|y_0\|^2_{L^2(\om)},\eqno{(2.4)}
$$
such that $P_r(y(\cdot,T))=0$, where $y(x,t)$ is the solution of the
equation (2.1) with $\delta=0$ and $u=u_r$,  and  where $C_1$ and
$C_2$ are the positive constants given in Theorem 2.2.}

\medskip
{\bf Proof}: Let $y(x,t)$ be the solution of the equation (2.1) with
$\delta=0$ and let $\varphi(x,t)$ be a solution of the equation
(2.2). Then
$$
<y(\cdot,T),\;\varphi(\cdot,T)>-<y_0(\cdot),\;\varphi(\cdot,0)>
=\d\int_0^T\int_\om\xe(t)\xw(x) u(x,t) \varphi(x,t) dxdt.
$$
Here and in what follows, $<\cdot,\;\cdot>$ denotes the inner
product in $L^2(\om)$. If we can show that
$<y(\cdot,T),\;\varphi(\cdot,T)>=0$ for all $\varphi(x,T)\in {\bf
X_r}$, then $P_r(y(\cdot,T))=0$. Thus, it suffices to prove that
there exists a control $u_r\in L^\infty(0,T;L^2(\om))$ with the
estimate (2.4) such that
$$
-<y_0(\cdot),\;\varphi(\cdot,0)>=\d\int_0^T\int_\om
\xe(t)\xw(x)u_r(x,t)\varphi(x,t)dxdt\;\;\mbox{for
all}\;\varphi(\cdot,T)\in {\bf X_r}.
$$
Now, we set
$$
{\bf Y_r}=\{\xe(t)\xw(x)\varphi(x,t);\;\varphi(x,t)\;\mbox{is the
solution to the equation (2.2) with }\varphi(\cdot,T)\in {\bf
X_r}\}.
$$
It is clear that ${\bf Y_r}$ is a linear subspace of
$L^1(0,T;L^2(\om))$. We define a linear functional $F_r:\;{\bf Y_r
}\ra {\r^1}$ by
$F_r(\xe\xw\varphi)=-<y_0(\cdot),\;\varphi(\cdot,0)>$.
 By the inequality (2.3), we see that
$$
\begin{array}{ll}
|F_r(\xe\xw\varphi)|^2&\leq
\|y_0\|^2_{L^2(\om)}\cdot\|\varphi(\cdot,0)\|^2_{L^2(\om)}\\
&\leq \d\f{C_1 e^{C_2\sqrt r}}{(m(E))^2}\|y_0\|^2_{L^2(\om)}\cdot
\|\xe\xw\varphi\|^2_{L^1(0,T;L^2(\om))}.
\end{array}
$$
Namely,
$$
\|F_r\|^2\leq \d\f{C_1 e^{C_2\sqrt
r}}{(m(E))^2}\|y_0\|^2_{L^2(\om)},
$$
where $\|F_r\|$ denotes the operator norm of $F_r$. Thus, $F_r$ is a
bounded linear functional on ${\bf Y_r}$. By the Hahn-Banach
Theorem, there is a bounded linear functional
$$
G_r:\;L^1(0,T;L^2(\om))\ra \r^1
$$
 such that
 $$
 G_r=F_r\mbox{ on }{\bf Y_r},
$$
and such that
$$
\|G_r\|^2=\|F_r\|^2\leq \d\f{C_1 e^{C_2\sqrt
r}}{(m(E))^2}\|y_0\|^2_{L^2(\om)}.
$$

Then, by making use of the Riesz Representation Theorem in [2], (
See p.61, [2].) there exists a function $u_r$ in the space
$L^\infty(0,T;L^2(\om))$ such that
$$
G_r(f)=\d\int_0^T\int_\om fu_rdxdt\;\;\;\;\mbox{for all }f\in
L^1(0,T;L^2(\om)),
$$
and such that
$$
\|u_r\|^2_{L^\infty(0,T;L^2(\om))}=\|G_r\|^2\leq \d\f{C_1
e^{C_2\sqrt r}}{(m(E))^2}\|y_0\|^2_{L^2(\om)}.
$$
In particular,
$$
F_r(\xe\xw\varphi)=\d\int_0^T\int_\om \xe\xw\varphi u_r
dxdt\;\;\mbox{for all }\xe\xw\varphi\in{\bf Y_r}.
$$
Namely,
$$
-<y_0(\cdot),\;\varphi(\cdot,0)>=\d\int_0^T\int_\om \xe\xw\varphi
u_rdxdt\;\;\mbox{for all }\varphi(\cdot,T)\in {\bf X_r}.
$$
This completes the proof.\endpf

\bigskip
The following lemma from the measure theory will be used in our
later discussion, whose  proof can be found in [11]. ( See p.
256-257, [11].)

\medskip
{\bf Lemma 2.4}. {\it For almost all $\w t$ in the set $E$, there
exists a sequence of numbers $\{t_i\}^\infty_{i=1}$ in the interval
$[0,T]$ such that
$$
t_1<\cdots<t_i<t_{i+1}<\cdots<\w t,\;t_i\ra \w t\;\mbox{as }i\ra
\infty,\eqno{(2.5)}
$$
$$
m(E\cap [t_i,t_{i+1}])\geq
\rho(t_{i+1}-t_i),\;i=1,2,\cdots,\eqno{(2.6)}
$$
and
$$
\d\f{t_{i+1}-t_i}{t_{i+2}-t_{i+1}}\leq
C_0,\;i=1,2,\cdots,\eqno{(2.7)}
$$
where $\rho$ and $C_0$ are two positive constants.}

\medskip

Now we are going to prove Theorem 2.1. Before proceeding the proof,
we  introduce briefly our main  strategy.  By applying Lemma 2.4,
there exist  a number $\w t$ and   a sequence $\{t_N\}^\infty_{N=1}$
 in the interval  $(0,T)$ such that (2.5)-(2.7) hold. The main part of the proof
is to show that for each $\tilde{y}_0$ in $L^2(\Omega)$, there
exists a control $\tilde{u}$ in the space $L^\infty ( t_1,
\tilde{t}; L^2(\Omega))$ with the estimate
$\|\tilde{u}\|^2_{L^\infty ( t_1, \tilde{t}; L^2(\Omega))}\leq L
\|\tilde{y}_0\|^2_{L^2(\Omega)}$ for a certain positive constant $L$
independent of $\tilde{y}_0$,  such that the solution
$\tilde{y}(x,t)$ to the equation:
$$
\left\{\begin{array}{ll} {\w y}_t(x,t) -\Delta \w
y(x,t)=\xe(t)\xw(x)\w u(x,t)\;\;\;\;\;\;\;\;\;\;\;\;\;\;\;&
\mbox{in }\;\;\om\times(\;t_1,\w t \;),\\
\w y(x,t)=0 & \mbox{on }\;\; \p\om\times(\;t_1,\w t\;),\\
\w y(x,t_1)=\w y_0(x)&\mbox{in }\;\;\om,
\end{array}\right.\eqno{(2.8)}
$$
has zero value at time $\tilde{t}$, namely,
$\tilde{y}(x,\tilde{t})=0$ over $\Omega$. To this end, we write
$$
[t_1,\tilde{t})=\bigcup^\infty_{N=1}(I_N\cup J_N),
$$
where $I_N=[t_{2N-1},t_{2N}]$ and
$J_N=[t_{2N},t_{2N+1}],\;N=1,2,\cdots$. Then we choose a suitable
sequence of positive numbers $\{r_N\}^\infty_{N=1}$ having the
following properties:\\
\hspace*{0.5 cm}\noindent (${\bf a}$) $ r_1< r_2<\cdots<r_N<\cdots,$
\\
\hspace*{0.5 cm}\noindent (${\bf b}$) $r_N\ra \infty$ as $N\ra
\infty$.\\
 On the subinterval $I_N$, we
control the heat equation  with a control $u_N$ restricted on the
subdomain $\omega\times(I_N\cap E)$ such that
$P_{r_N}(y_N(\cdot,t_{2N}))=0$, where $P_{r_N}$ denotes the
orthogonal projection from $L^2(\om)$ onto
$\mbox{span}\;\{X_i(x)\}_{i=1}^{r_N}$. On the subinterval $J_N$, we
let the heat equation to evolve freely.  We start with the initial
data for the equation on $I_1$ to be $y_0$. For the  initial data on
$I_N$, $ N=2,3,\cdots$, we define it  to be the ending value of the
solution for the equation on $J_{N-1}$. The  initial data of the
equation on $J_N$, $N=1,2, \cdots$, is given by the ending value of
the solution for the equation on $I_N$.  Moreover, by making use of
Lemma 2.3 and Lemma 2.4, we will show that there is a sequence
$\{r_N\}^\infty_{N=1}$, having the properties $\bf (a)$ and $\bf
(b)$ as above,  such that the $L^\infty(I_N;L^2(\om))$-norm of the
control $u_N$ is bounded by
$L^{\f{1}{2}}\|\tilde{y}_0\|_{L^2(\Omega)}$ for a certain positive
constant $L$ independent of $N$ and $\tilde{y}_0$.  Then, we
construct a control $\tilde{u}$ by setting
$$
 \tilde{u}(x,t)=\left \{\begin{array}{ll}
u_N(x,t),\;\;\;\;&x\in\om,\;t\in I_N,\;N=1,2,\cdots,\\
0,&x\in \om ,\;t\in J_N,\;N=1,2,\cdots.
\end{array}\right.
$$
We can show that  this control $\tilde{u}$ makes the corresponding
trajectory $\tilde{y}$ of the equation (2.8) have zero value at time
$\tilde{t}$.

Now, we set
$$
 u(x,t)=\left \{\begin{array}{ll}
\tilde{u}(x,t),\;\;\;\;&   \mbox{in}\; \om\times (t_1, \tilde{t} ),\\
0,&      \mbox{in}\; \om\times ( (0,T)\backslash (t_1,\tilde{t} ) )
\end{array}\right.
$$
and take $\tilde{y}_0$ to be $\psi(x, t_1)$, where $\psi(x,t)$ is
the solution of the heat equation on $\Omega\times (0, t_1)$ with
the initial data $y_0$. Then it is clear that this control $u$ makes
the trajectory $y(x,t)$ of the equation (2.1) with $\delta=0$ have
zero value at time $T$. Moreover, $\|u\|^2_{L^\infty(0,T;
L^2(\Omega))}\leq L\|y_0\|^2_{L^2(\Omega)}$.

 We next replace the
sequence $\{t_N\}_{N=1}^\infty$ and the number $\tilde{t}$ by the
sequence $\{t_N-\delta\}_{N=1}^\infty$ and the number
$(\tilde{t}-\delta)$ respectively, where the number $\delta$ is such
that $0\leq\delta\leq t_1$. Then by making use of  the same argument
as above, we obtain that for each number $\delta$ with
$0\leq\delta\leq t_1$, there exists a control $u_\delta$ in the
space $L^\infty(0, T-\delta; L^2(\Omega))$ with the estimate
$\|u_\delta\|^2_{L^\infty(0,T-\delta; L^2(\Omega))}\leq
L_\delta\|y_0\|^2_{L^2(\Omega)}$ for a certain positive number
$L_\delta$ independent of $y_0$, such that the corresponding
solution $y^\delta$ to the equation (2.1) reaches zero value at time
$T-\delta$, namely, $y^\delta(x, T-\delta)=0$ over $\Omega$. We
finally prove that $L_\delta=L$ is independent of $\delta$.

Now we turn to prove Theorem 2.1.

\medskip

{\bf Proof of Theorem 2.1.}  Without loss of generality, we can
assume that $C_1\geq 1$, where $C_1$ is the positive constant given
in Theorem 2.2. By making use of Lemma 2.4, we can take a number $\w
t$ in the set $E$ with $\w t<T$ and a sequence
$\{t_N\}^\infty_{N=1}$ in the open interval $(0,T)$ such that
(2.5)-(2.7) hold for certain positive numbers $\rho$ and $C_0$  and
such that
$$
\w t -t_1\leq \;\mbox{Min}\{\l_1,\;1\}.
$$
We shall first prove that for each $\w{y}_0$ in $L^2(\om)$, there
exists a control  $\w u$ in the space $L^\infty(t_1,\w t;L^2(\om))$
with the estimate $ \|\w u\|^2_{L^\infty(t_1,\w t;L^2(\om))}\leq L
\|\w y_0\|^2_{L^2(\om)}$ for a certain positive constant $L$
independent of $\w{y}_0$, such that the solution $\w y$ to the
equation (2.8)
reaches zero value at time $\w t$, namely,  $\w y (x,\w t \;)=0$ over $\om$.\\
\hspace*{0.5 cm} To this end, we shall use the  strategy presented
above. We set $I_N=[t_{2N-1},t_{2N}]$, $J_N=[t_{2N},t_{2N+1}]$ for
$N=1,2,\cdots$. Then
$$
[t_1,\w t\;)=\bigcup^\infty_{N=1}(I_N\bigcup J_N).
$$
Notice  that for each $N\geq 1$,  it holds that $m(E\cap I_N)>0$.

Now, on the interval $I_1\equiv [t_1,t_2]$, we consider the
following controlled heat equation:
$$
\left\{\begin{array}{ll} y'_1(x,t) -\Delta y_1(x,t) =\xe(t)\xw(x)
u_1(x,t) \;\;\;\;\;\;\;\;\;\;\;\;\;\;\;&
\mbox{in }\;\;\om\times(\;t_1,t_2 \;),\\
y_1(x,t)=0 & \mbox{on }\;\; \p\om\times(\;t_1,t_2\;),\\
y_1(x,t_1)=\w y_0(x)&\mbox{in }\;\;\om.
\end{array}\right.
$$
By Lemma 2.3, for any $r_1>0$, there exists a control $u_1$ in the
space $L^\infty (t_1,t_2;L^2(\om))$ with the estimate:
$$
\begin{array}{ll}
\|u_1\|^2_{L^\infty(t_1,t_2;L^2(\om))}&\leq \d\f{C_1 e^{C_2\sqrt
r_1}}{(m(E\cap [t_1,t_2]))^2}\|\w y_0\|^2_{L^2(\om)},\\
\end{array}
$$
such that $P_{r_1}(y_1(\cdot,t_2))=0$. Then, by (2.6) and (2.7) in
Lemma 2.4, we see that
$$
\begin{array}{ll}
\|u_1\|^2_{L^\infty(t_1,t_2;L^2(\om))}&\leq  \d\f{C_1 e^{C_2\sqrt
r_1}}{\rho^2(t_2-t_1)^2}\|\w
y_0\|^2_{L^2(\om)}\\
&\equiv\d\f{C_1}{\rho^2(t_2-t_1)^2}\cdot \alpha_1 \|\w
y_0\|^2_{L^2(\om)},
\end{array}
$$
where $\alpha_1=e^{C_2\sqrt r_1}$. Moreover, we have
$$
\begin{array}{ll}
\|y_1(\cdot,t_2)\|^2_{L^2(\om)}&\leq
\|y_1(\cdot,t_1)\|^2_{L^2(\om)}
+\f{1}{\l_1}\d\int^{t_2}_{t_1}\|u_1(\cdot,s)\|^2_{L^2(\om)}ds\\
&\leq\|\w y_0\|^2_{L^2(\om)}
+\f{(t_2-t_1)}{\l_1}\|u_1\|^2_{L^\infty(t_1,t_2;L^2(\om))}\\
&\leq 2\d\f{C_1}{\rho^2(t_2-t_1)^2}\cdot \alpha_1\|\w
y_0\|^2_{L^2(\om)}.
\end{array}
$$
Here we have used the facts that $(t_2-t_1)\leq
\mbox{Min}\;(\l_1,1),\;\rho<1$ and $C_1>1$.

On the interval $J_1 \equiv [t_2,t_3]$, we consider the following
heat equation without control:
$$
\left\{\begin{array}{ll} z'_1(x,t) -\Delta z_1(x,t) =0
\;\;\;\;\;\;\;\;\;\;\;\;\;\;\;&
\mbox{in }\;\;\om\times(t_2,t_3 ),\\
z_1(x,t)=0 & \mbox{on }\;\; \p\om\times(t_2,t_3),\\
z_1(x,t_2)=y_1(x,t_2)&\mbox{in }\;\;\om.
\end{array}\right.
$$
Since $P_{r_1}(y_1(\cdot,t_2))=0$, we have
$$
\begin{array}{ll}
\|z_1(\cdot,t_3)\|^2_{L^2(\om)}&\leq
\mbox{exp}\;(-2r_1(t_3-t_2))\cdot\|y_1(\cdot,t_2)\|^2_{L^2(\om)}\\
&\leq 2\d\f{C_1}{\rho^2(t_2-t_1)^2}\alpha_1\cdot
\mbox{exp}\;(-2r_1(t_3-t_2))\cdot\|\w y_0\|^2_{L^2(\om)}.
\end{array}
$$

On the interval $I_2\equiv [t_3,t_4]$, we consider the controlled
heat equation as follows:
$$
\left\{\begin{array}{ll} y'_2(x,t) -\Delta y_2(x,t) =\xe(t)\xw(x)
u_2(x,t) \;\;\;\;\;\;\;\;\;\;\;\;\;\;\;&
\mbox{in }\;\;\om\times(t_3,t_4),\\
y_2(x,t)=0 & \mbox{on }\;\; \p\om\times(t_3,t_4),\\
y_2(x,t_3)=z_1(x,t_3)&\mbox{in }\;\;\om.
\end{array}\right.
$$
Then by Lemma 2.3, for any $r_2>0$, there exists a control $u_2$ in
the space $L^\infty (t_3,t_4;L^2(\om))$ with the estimate:
$$
\|u_2\|^2_{L^\infty(t_3,t_4;L^2(\om))}\leq \d\f{C_1 e^{C_2\sqrt
{r_2}}}{m(E\cap [t_3, t_4]))^2}\cdot\|z_1(\cdot,t_3)\|^2_{L^2(\om)},
$$
such that $P_{r_2}(y_2(\cdot,t_4))=0$. By (2.6) and (2.7) in Lemma
2.4, we get
$$
\|u_2\|^2_{L^\infty(t_3,t_4;L^2(\om))}\leq
2(\d\f{C_1}{\rho^2(t_2-t_1)^2})^2
C_0^4\cdot\alpha_1\cdot\alpha_2\cdot\|\w y_0\|^2_{L^2(\om)}
$$
where $\alpha_2=\mbox{exp }(C_2\sqrt {r_2})\mbox{exp
}(-2r_1(t_3-t_2))$. Moreover, it holds that
$$
\begin{array}{ll}
\|y_2(\cdot,t_4)\|^2_{L^2(\om)}&\leq \|z_1(\cdot,t_3)\|^2_{L^2(\om)}
+\d\f{1}{\l_1}(t_4-t_3)\|u_2\|^2_{L^\infty(t_3,t_4;L^2(\om))}\\
&\leq 2^2(\d\f{C_1}{\rho^2(t_2-t_1)^2})^2
C_0^4\cdot\alpha_1\cdot\alpha_2\cdot\|\w y_0\|^2_{L^2(\om)}.
\end{array}
$$

On the interval $J_2\equiv [t_4,t_5]$, we consider the following
heat equation without control:
$$
\left\{\begin{array}{ll} z'_2(x,t) -\Delta z_2(x,t) =0
\;\;\;\;\;\;\;\;\;\;\;\;\;\;\;&
\mbox{in }\;\;\om\times(t_4,t_5 ),\\
z_2(x,t)=0 & \mbox{on }\;\; \p\om\times(t_4,t_5),\\
z_2(x,t_4)=y_2(x,t_4)&\mbox{in }\;\;\om.
\end{array}\right.
$$
Since $P_{r_2}(y_2(\cdot,t_4))=0$, we have
$$
\begin{array}{ll}
\|z_2(\cdot,t_5)\|^2_{L^2(\om)}&\leq \mbox{exp
}(-2r_2(t_5-t_4))\|y_2(\cdot,t_4)\|^2_{L^2(\om)}\\
&\leq 2^2(\d\f{C_1}{\rho^2(t_2-t_1)^2})^2
C_0^4\cdot\alpha_1\cdot\alpha_2\cdot\|\w
y_0\|^2_{L^2(\om)}\cdot\mbox{exp }(-2r_2(t_5-t_4)).
\end{array}
$$

On  the interval $I_3\equiv [t_5,t_6]$, we consider the following
controlled heat equation:
$$
\left\{\begin{array}{ll} y'_3(x,t) -\Delta y_3(x,t) =\xe(t)\xw(x)
u_3(x,t) \;\;\;\;\;\;\;\;\;\;\;\;\;\;\;&
\mbox{in }\;\;\om\times(t_5,t_6 ),\\
y_3(x,t)=0 & \mbox{on }\;\; \p\om\times(t_5,t_6),\\
y_3(x,t_5)=z_2(x,t_5)&\mbox{in }\;\;\om.
\end{array}\right.
$$
Then by Lemma 2.3, for any $r_3>0$, there exists a control $u_3$ in
the space $L^\infty (t_5,t_6;L^2(\om))$ with the estimate:
$$
\|u_3\|^2_{L^\infty(t_5,t_6;L^2(\om))}\leq \d\f{C_1 e^{C_2\sqrt
{r_3}}}{(m(E\cap [t_5,t_6]))^2}\|z_2(\cdot,t_5)\|^2_{L^2(\om)},
$$
such that $P_{r_3}(y_3(\cdot,t_6))=0$. By making use of (2.6) and
(2.7) again, we get
$$
\|u_3\|^2_{L^\infty(t_5,t_6;L^2(\om))} \leq
2^2(\d\f{C_1}{\rho^2(t_2-t_1)^2})^3 C_0^4\cdot C_0^{4\cdot
2}\cdot\alpha_1\cdot\alpha_2\cdot\alpha_3\cdot\|\w
y_0\|^2_{L^2(\om)},
$$
where $\alpha_3=\mbox{exp }(C_2\sqrt {r_3})\;\mbox{exp
}(-2r_2(t_3-t_2)C_0^{-2})$.

\medskip

Generally, on the interval $I_N$, we consider the controlled heat
equation:
$$
\left\{\begin{array}{ll} y'_N(x,t)-\Delta y_N(x,t) =\xe(t)\xw(x)
u_N(x,t) \;\;\;\;\;\;\;\;\;\;\;\;\;\;\;&
\mbox{in }\;\;\om\times(t_{2N-1},t_{2N}),\\
y_N(x,t)=0 & \mbox{on }\;\; \p\om\times(t_{2N-1},t_{2N}),\\
y_N(x,t_{2N-1})=z_{N-1}(x,t_{2N-1})&\mbox{in }\;\;\om.
\end{array}\right.
$$
On the interval $J_N$, we consider the following heat equation
without control:
$$
\left\{\begin{array}{ll} z'_N(x,t) -\Delta z_N(x,t) =0
\;\;\;\;\;\;\;\;\;\;\;\;\;\;\;&
\mbox{in }\;\;\om\times(t_{2N},t_{2N+1} ),\\
z_N(x,t)=0 & \mbox{on }\;\; \p\om\times(t_{2N},t_{2N+1}),\\
z_N(x,t_{2N})=y_N(x,t_{2N})&\mbox{in }\;\;\om.
\end{array}\right.
$$
Then by making use of induction argument, we can obtain the
following: {\it For each $r_N>0$, there exists a control $u_N$ in
the space $L^\infty(I_N;L^2(\om))$ with the following estimate:
$$
\begin{array}{ll}
\|u_N\|^2_{L^\infty(I_N;L^2(\om))}\\
\leq 2^{N-1}(\d\f{C_1}{\rho^2(t_2-t_1)^2})^N C_0^4\cdot C_0^{4\cdot
2}\cdots C^{4(N-1)}_0\cdot
\alpha_1\cdot\alpha_2\cdots\alpha_N\cdot\|\w y_0\|^2_{L^2(\om)},
\end{array}
$$
where
$$
\alpha_N=\left\{\begin{array}{ll} \mbox{exp }(C_2 \sqrt
r_1),\;\;\;\;\;\;\;\;\;\;\;\;\;&N=1,\\
 \mbox{exp }(C_2\sqrt
{r_N})\;\mbox{exp }(-2r_{N-1}(t_3-t_2)C_0^{-2(N-2)}),&N\geq 2,
\end{array}\right.
\eqno{(2.9)}
$$
such that $P_{r_N}(y_N(\cdot,t_{2N}))=0$.}  It is easily seen that
for each $N\geq 1$,
$$
\|u_N\|^2_{L^\infty(I_N;L^2(\om))}\leq(\w
C)^{N(N-1)}\alpha_1\cdots\alpha_N\cdot\|\w y_0\|^2_{L^2(\om)},
\eqno{(2.10)}
$$
where
$$
\w C=\d\f{2C_1}{\rho^2(t_2-t_1)^2}\cdot C^2_0.
 \eqno{(2.11)}
$$

Now, we set
$$ r_{N}=[\d\f{2}{(t_3-t_2)}{\w C}^{N-1}]^4\equiv [
A\cdot {\w C}^{N-1}]^4, \;\; N\geq 1.
\eqno{(2.12)}
$$
Because we have $\w C>C^2_0>1$ and $t_3-t_2<1$, it holds that
$$
2^4<r_1<r_2<\cdots <r_N<r_{N+1}<\cdots,\;\; \mbox{and}\;\;
r_N\ra\infty\;\;\mbox{as}\; N\ra\infty.
$$
Moreover, we have
$$
{r_{N-1}}^{\f{1}{4}}(t_3-t_2)C_0^{-2(N-2)}\geq 2\;\;\; \mbox{for
each}\; N\geq 2.
$$
Then we get
$$
\mbox{exp }\{-2r_{N-1}(t_3-t_2)C_0^{-2(N-2)}\}\leq \mbox{exp
}(-4{r_{N-1}}^{\f{3}{4}})\;\;\mbox{for each}\; N\geq 2.
\eqno{(2.13)}
$$
Since
$$
\begin{array}{ll}
{\w C}^{N(N-1)}\mbox{exp }(-{r_{N-1}}^{\f{3}{4}})=\d\f{{\w
C}^{N(N-1)}}
{(\mbox{exp}({r_{N-1}}^{\f{1}{4}}))^{{r_{N-1}}^{\f{1}{2}}} } &\leq
\d\f{{\w C}^{N(N-1)}}{(\mbox{exp }(2{\w
C}^{N-1}))^{{r_{N-1}}^{\f{1}{2}}}}\\
&\leq \d\f{{\w C}^{N(N-1)}}{{\w C}^{(N-1)\cdot 2\cdot
{r_{N-1}}^{\f{1}{2}}}}
\end{array}
$$
 for each $N\geq 2$, we derive from (2.12) that there exists a natural number $N_1$ with
$N_1\geq 2$ such that for each $N\geq N_1$,
$$
{\w C}^{N(N-1)}\mbox{exp }(-{r_{N-1}}^{\f{3}{4}})\leq 1.
\eqno{(2.14)}
$$
By making use of  (2.12) again, we obtain that for each $N\geq 2$,
$$
\mbox{exp }(C_2\sqrt r_N)\;\mbox{exp
}(-{r_{N-1}}^{\f{3}{4}})=\mbox{exp }(C_2 A^2 \w
C^{2(N-1)})\;\mbox{exp }(-A^3 \w C^{3(N-2)}).
$$
Thus, there exists a natural number $N_2$ with $N_2\geq 2$ such that
for each $N\geq N_2 $,
$$
\mbox{exp }(C_2\sqrt r_N)\;\mbox{exp }(-{r_{N-1}}^{\f{3}{4}})\leq
1.\eqno{(2.15)}
$$
Now we set
$$
N_0=\mbox{max }\{N_1,N_2\}.
 \eqno{(2.16)}
$$
 Then by (2.13), (2.14) and
(2.15), we see that for all $N\geq N_0$,
$$
\begin{array}{ll}
\w C^{N(N-1)}\alpha_N\\
=\w C^{N(N-1)}\mbox{exp }(C_2\sqrt r_N)\;\mbox{exp
}(-2r_{N-1}(t_3-t_2)C_0^{-2(N-2)})\\\leq \w C^{N(N-1)}\mbox{exp
}(C_2\sqrt r_N)\mbox{exp }(-4{r_{N-1}}^{\f{3}{4}})\\
\leq \mbox{exp }(-2{r_{N-1}}^{\f{3}{4}}).
\end{array}\eqno{(2.17)}
$$
Moreover, it is obvious that
$$
\alpha_N\leq 1\;\;\mbox{for all }N\geq N_0. \eqno{(2.18)}
$$
Now, we set
$$
L=\mbox{max}\; \{\;(\w C)^{N(N-1)}\alpha_1\cdots\alpha_N, \;\; 1\leq
N\leq N_0 \;\}. \eqno{(2.19)}
$$
It follows from (2.10), (2.17), (2.18) and (2.19) that  for all
$N\geq 1$,
$$
\|u_N\|^2_{L^\infty(I_N;L^2(\om))}\leq L \|\w
y_0\|^2_{L^2(\om)}.\eqno{(2.20)}
$$
Then we construct a control $\w u$ by setting
$$ \w
u(x,t)=\left\{\begin{array}{ll} u_N(x,t),\;\;\;\;&x\in \om,\;t\in
I_N,\;N\geq 1,\\
0,&x\in \om,\;t\in J_N,\;N\geq 1,
\end{array}\right.
\eqno{(2.21)}
$$
from which and by (2.20), we easily see that the control $\w u$ is
in the space $L^\infty (t_1,\w t;L^2(\om))$ and satisfies the
estimate:
$$
\|\w u\|^2_{L^\infty(t_1,\w t;L^2(\om))}\leq L \|\w
y_0\|^2_{L^2(\om)}.
$$
Let $\w y$ be the solution of the  equation (2.8) corresponding to
the control $\w u$ constructed in (2.21). Then on the interval
$I_N$, $\w y(\cdot,t)=y_N(\cdot,t)$. Since
$P_{r_N}(y_N(\cdot,t_{2N}))=0$ for all $N\geq 1$ and
$r_1<r_2<\cdots<r_N<\cdots$,  by making use of (2.21) again, we see
that
$$
P_{r_N}(\w y (\cdot , t_{2M}))=0\;\;\mbox{for all}\;\;  M\geq N.
\eqno{(2.22)}
$$
On the other hand, since $t_{2M}\ra \w t$ as $ M\ra\infty$, we
obtain that
$$
\w y(\cdot, t_{2M})\ra \w y(\cdot, \w t\; )\;\;\mbox{strongly in}\;
L^2(\Omega),\;\;\mbox{as}\; M\ra\infty.
$$
This, together with (2.22), implies that $P_{r_N}( \w y(\cdot, \w
t\; ) )=0$ for all $N\geq 1$. Since $r_N\ra \infty$ when
$N\ra\infty$, it holds that  $\w y(\cdot, \w t \;)=0$. Thus, we have
proved that for each $\w{y}_0\in L^2(\om)$, there exists a control
$\w u\in L^\infty(t_1,\w t;L^2(\om))$ with the estimate $ \|\w
u\|^2_{L^\infty(t_1,\w t;L^2(\om))}\leq L \|\w y_0\|^2_{L^2(\om)}$,
where the constant $L$ is given by (2.19), such that the solution
$\w y$ to the equation (2.8) reaches zero value at time $\w t$,
namely, $\w y (x,\w t \;)=0$ over $\om$.

Now, we take $\w y_0(x)$ to be $\psi(x,t_1)$, where $\psi(x,t)$ is
the solution to the following equation:
$$
\left\{\begin{array}{ll}
\psi_t(x,t)-\Delta\psi(x,t)=0\;\;\;\;\;\;\;\;&\mbox{in } \om\times
(0,t_1),\\
\psi(x,t)=0&\mbox{on }\p\om\times(0,t_1),\\
\psi(x,0)=y_0(x)  &\mbox{in }\om
\end{array}\right.
$$
and construct a control $u$ by setting
$$
u(x,t)=\left\{\begin{array}{ll} 0\;\;\;\;\;\;\;\;\;\;\;\;&\mbox{in }
\om\times
(0,t_1),\\
\w u(x,t)&\mbox{in }\om\times(t_1,\w t\;),\\
0 &\mbox{in }\om\times(\w t,T).
\end{array}\right.
\eqno{(2.23)}
$$
It is clear that this control  $u$ is in the space
$L^\infty(0,T;L^2(\om))$ and that the  corresponding  solution $y$
of  the equation (2.1) with $\delta=0$ reaches zero value at time
$T$, namely, $y(x,T)=0$ over $\om$. Moreover, the control $u$
constructed in (2.23) satisfies the following estimate:
$$
\| u\|^2_{L^\infty(0,T;L^2(\om))}\leq L \| y_0\|^2_{L^2(\om)},
$$
where $L$ is given by (2.19).

Next, we take $\delta_0$ to be the number $t_1$ given above.
   For each  $\delta$ with
 $0\leq \delta\leq\delta_0$, we set
$$
\w {t_\delta}= \w t-\delta\;\; \mbox{and}\;\;
t_{N,\delta}=t_N-\delta\;\;\mbox{for all}\;\; N=1,2, \cdots.
$$
Then it holds that
$$
 0\leq t_{1,\delta}< t_{2,\delta}<\cdots <t_{N,\delta}\ra \w{t_\delta}< T-\delta.
 $$
Moreover, we have  for each $N\geq 1$,
$$
 m(E_\delta \cap [t_{N,\delta}, t_{N+1, \delta}])=m(E\cap [t_N,t_{N+1}])\geq
\rho(t_{N+1}-t_N),
$$
and
$$
\d\f{t_{N+1, \delta}-t_{N,\delta}}{t_{N+2,
\delta}-t_{N+1,\delta}}=\d\f{t_{N+1}-t_{N}}{t_{N+2}-t_{N+1}}
  \leq C_0,
$$
where $C_0$ and $\rho$ are the positive constants as above.

Now, we can use  exactly the same argument as above to get for each
$\delta$ with $0\leq\delta\leq\delta_0$, the existence of a control
$u_\delta(t)$ in the space $L^\infty (0, T-\delta; L^2(\Omega))$
such that the corresponding solution $y^\delta$ to the equation
(2.1) reaches zero value at time $T-\delta$, namely, $y^\delta(x,
T-\delta)=0$ over $\Omega$. Moreover, this control $u_\delta$
satisfies the following estimate: ( See (2.9)-(2.12) and (2.19).)
$$
\| u_\delta\|^2_{L^\infty(0,T-\delta ;L^2(\om))}\leq L_\delta\cdot
\| y_0\|^2_{L^2(\om)}.
$$
The constant $L_\delta$ is given by
$$
L_\delta=\mbox{max}\; \{({\w C}_\delta)^{N(N-1)}
\alpha_{1,\delta}\cdots\alpha_{N,\delta}, \;\; 1\leq N\leq N_0 \},
$$
where
$$
{\w C}_\delta=\d\f{2C_1}{\rho^2(t_{2,\delta}-t_{1,\delta})^2}\cdot
C^2_0
 $$
and
$$
\alpha_{N,\delta}=\left\{\begin{array}{ll} \mbox{exp }(C_2 \sqrt
{r_{1,\delta}}),\;\;\;\;\;\;\;\;\;\;\;\;\;&N=1,\\
 \mbox{exp }(C_2\sqrt
{r_{N,\delta}})\;\mbox{exp
}(-2r_{N-1,\delta}(t_{3,\delta}-t_{2,\delta})C_0^{-2(N-2)}),&N\geq
2,
\end{array}\right.
$$
with
$$
r_{N,\delta}=[\d\f{2}{(t_{3,\delta}-t_{2,\delta})}{\w
C}_\delta^{N-1}]^4,\;\; N=1,2,\cdots,
$$
and where the natural number $N_0$ is  given by (2.16). Since
$$
t_{N+1,\delta}-t_{N,\delta}= t_{N+1}-t_N, \;\;\mbox{for all}\;\;
N=1,2,\cdots,
$$
we see easily that $\w C_\delta =\w C$ and
$\alpha_{N,\delta}=\alpha_N$ for all $N\geq 1$. Then  it holds that
$ L_\delta = L$ for all $\delta$ with $0\leq \delta \leq \delta_0$.
This completes the proof.\endpf


\section{The bang-bang principle for time optimal control}

\hspace*{0.5 cm} In this section, we shall prove the main result of
the paper, namely,  each optimal control for the problem $\bf (P)$
satisfies the bang-bang principle in the weaker form. Moreover, we
shall show the uniqueness of the optimal control for the problem
$\bf (P)$, when the target set $S$ is  convex  and the control set
is a closed ball. Throughout of this section, we shall denote by
$y(t;u, y_0)$ the solution of the equation (1.1) corresponding to
the control $u$ and the initial data $y_0$,  and write
$\{G(t)\}_{t\geq 0}$ for the semigroup generated by $\Delta$ with
the Dirichlet boundary condition.

\medskip
{\bf Theorem 3.1}. {\it Suppose that the control set $U$ is closed,
bounded and nonempty  in $L^2(\Omega)$ and the target set $S$ is
 nonempty in $L^2(\Omega)$. Let $T^*$ be the optimal time
 and   $u^*$ be an optimal
control for  the problem $\bf (P)$. Then it holds that $u^*(t)\in \p
U$ for almost all $t\in[0,T^*]$. If we further assume that
$\chi_{\omega} U\subset U$, then it holds that
$\chi_{\omega}u^*(t)\in \p U$ for almost all $t\in [0, T^*]$.}

\medskip
{\bf Proof of Theorem 3.1}. Seeking a contradiction, we suppose that
there exist a  subset $E$ of positive measure in the interval
$[0,T^*]$  and a positive number $\e$ such that the following holds:
$$
u^*(t)\in U\;\;\mbox{and}\; \;\mbox{d} ( u^*(t), \p U)\geq \e\;\;
\mbox{for each}\; t \;\mbox{in the set}\; E,
$$
where $\mbox{d} (u^*(t), \p U)$ denotes the distance of the point
$u^*(t)$ to the set $\p U$ in $L^2(\Omega)$. Then we  would get
$$
B( u^*(t), \f{\e}{2}) \subset U\;\; \mbox{for each}\; t \;\mbox{in
the set}\;  E. \eqno{(3.1)}
$$

 We shall obtain from (3.1) that {\it there exist a positive number $\delta$
with $\delta<T^*$ and a control $v_\delta$ in the set $\uad$ such
that the following holds:
$$
y(T^*-\delta; v_\delta, y_0)=y(T^*; u^*,y_0). \eqno{(3.2)}
$$}
Thus, $T^*$ could not be the optimal time for the problem $\bf (P)$,
which leads to a contradiction.

We  first observe that
$$
\begin{array}{ll}
y(T^*-\delta; v_\delta, y_0)=G(T^*-\delta)y_0+\d\int_0^{T^*-\delta}
 G(T^*-\delta-\sigma)\xw v_\delta(\sigma)d\sigma,\\
y(T^*;u^*, y_0)=G(T^*)y_0+\d\int_0^{T^*} G(T^*-\sigma)\xw
u^*(\sigma)d\sigma.
\end{array}
$$
Hence,  (3.2) is equivalent to the following: {\it There exist a
positive number $\delta$ with $\delta<T^*$ and a control $v_\delta$
in the set $\uad$
 such that the following holds:
$$
\d\int_0^{T^*-\delta} G(T^*-\delta-\sigma)\xw
v_\delta(\sigma)d\sigma=[G(T^*)-G(T^*-\delta)]y_0+\d\int_0^{T^*}G(T^*-\sigma)\xw
u^*(\sigma)d\sigma.\eqno{(3.3)}
$$}
Notice that for any positive number $\delta$ with $\delta<T^*$, we
have
$$
\begin{array}{ll}
\d\int_0^{T^*} G(T^*- \sigma)\xw u^*(\sigma)d\sigma\\
=\d\int_0^{\delta}G(T^*-\sigma)\xw
u^*(\sigma)d\sigma+\d\int^{T^*}_{\delta}G(T^*-\sigma)\xw
u^*(\sigma)d\sigma\\
=G(T^*-\delta)\d\int_0^{\delta}G(\delta-\sigma)\xw
u^*(\sigma)d\sigma+\d\int_0^{T^*-\delta}G(T^*-\delta-\sigma)\xw
u^*(\delta+\sigma)d\sigma
\end{array}
$$
and
$$
[G(T^*)-G(T^*-\delta)]y_0=G(T^*-\delta)[(G(\delta)-I)y_0].
$$
Therefore, (3.3) is equivalent to the following: {\it There exist a
positive number $\delta$ with $\delta<T^*$ and a control $v_\delta$
in the set $\uad$
 such that the following holds:
$$
\begin{array}{ll}
\d\int_0^{T^*-\delta}G(T^*-\delta-\sigma)\xw
v_\delta(\sigma)d\sigma\\=G(T^*-\delta)[\d\int_0^{\delta}G(\delta-\sigma)\xw
u^*(\sigma)d\sigma +(G(\delta)-I)y_0]\\
+\d\int_0^{T^*-\delta}G(T^*-\delta-\sigma)\xw
u^*(\sigma+\delta)d\sigma\\
 \equiv G(T^*-\delta)h_\delta + \d\int_0^{T^*-\delta}G(T^*-\delta-\sigma)\xw
u^*(\sigma+\delta)d\sigma,
\end{array}\eqno{(3.4)}
$$
where
$$
h_\delta=\d\int_0^{\delta}G(\delta-\sigma)\xw
u^*(\sigma)d\sigma+(G(\delta)-I)y_0.\eqno{(3.5)}
$$}

For each positive number $\delta$, we write   $E_\delta$ for the set
$\{ t; \; t+\delta\in E \}$ and denote by $\chi_{E_\delta}$  the
characteristic function of the set $E_\delta$. We first claim  the
following: {\it For each positive number $\delta$ sufficiently
small, there exists a control $u_\delta$ in the space
$L^\infty(0,\infty; L^2(\om))$ such that
$$
\|u_\delta(t)\|_{L^2(\Omega)}\leq  \f{\e}{2}\;\;\;\mbox{for almost
all }\; t\geq 0, \eqno{(3.6)}
$$
and such that
$$
y(T^*-\delta; \chi_{E_\delta}u_\delta, 0)=G(T^*-\delta)h_\delta.
\eqno{(3.7)}
$$}Recall that $y(t; \chi_{E_\delta}u_\delta, 0)$ is the solution of
the controlled heat equation (1.1) with $u$ and $y_0$ being replaced
by $\chi_{E_\delta}u_\delta$ and $0$ respectively, and that
$\varphi(t)\equiv G(t)h_\delta$ is the solution of the equation
(1.1) with $u$ and $y_0$ being replaced by $0$ and $h_\delta$
respectively. Then, what we claimed above is obviously equivalent to
the following: {\it For each positive number $\delta$ sufficiently
small, there exists a control $u_\delta$ with the estimate:
$$
\|u_\delta(t)\|_{L^2(\Omega)}\leq \f{\e}{2}\;\; \mbox{for almost
all}\;\; t\geq 0,
$$
 such that the following holds:
$$ z^\delta(T^*-\delta)=0,$$
where $z^\delta(t)$ is the solution to the following controlled heat
equation:
$$
\left\{\begin{array}{ll}
z^\delta_t(t)-\Delta z^\delta(t)=\xw \chi_{E_\delta}(t) u_\delta(t)\;\;\;\mbox{in }(0,T^*-\delta),\\
z^\delta(0)=-h_\delta.
\end{array}\right.
\eqno{(3.8)}
$$}However, by Theorem 2.1,  there exist positive numbers $\delta_0$
and $L$ such that for each $\delta$ with $0<\delta\leq\delta_0$,
there is  a control $u_\delta$ in the space $L^\infty(0,T^*-\delta;
L^2(\Omega) )$ with the estimate:
$$
\|u_\delta\|^2_{L^\infty(0,T^*-\delta; L^2(\om))}\leq
L\|h_\delta\|^2_{L^2(\om)}, \eqno{(3.9)}
$$
such that the following holds:
$$
z^\delta(T^*-\delta)=0. \eqno{(3.10)}
$$
 On the other hand, by (3.5), we can get a positive number $\w\delta$ such that
for each positive number $\delta$ with $\delta\leq \w\delta$, the
following holds:
$$
\|h_\delta\|^2_{L^2(\Omega)}\leq {(\f{\e}{2})}^2/L.
$$
This, together with (3.9), implies that for each positive number
$\delta$ with $\delta\leq \mbox{min}\{ \delta_0, \w\delta\}$, there
is a control $u_\delta$ with the estimate:
$$
\|u_\delta\|_{L^\infty(0,T^*-\delta; L^2(\om))}\leq \f{\e}{2},
\eqno{(3.11)}
$$
such that the corresponding solution $z^\delta$ to the equation
(3.8) satisfies (3.10).

 Next, we fix such a positive number $\delta$  and the corresponding control $u_\delta$  that
(3.10) and (3.11) hold. Then we extend the control $u_\delta(\cdot)$
by setting it to be zero on the interval $(T^*-\delta, \infty)$, and
still denote the extension by $u_\delta(\cdot)$. Clearly, this
extended control $u_\delta$ is in the space $L^\infty(0, \infty;
L^2(\om))$ and makes (3.6) and (3.7) hold. Thus, we have proved the
above mentioned claim.

Now, we take an element $u_0$  from the control set $U$ and
construct a control $v_\delta$ by setting
$$
v_\delta(t)=\left\{\begin{array}{ll}u^*(t+\delta) +
\chi_{E_\delta}(t)u_\delta(t),\;\;&\mbox{if}\; t\in [0, T^*-\delta],\\
u_0, \;\; \;&\mbox{if}\; t>T^*-\delta.
\end{array}\right.\eqno{(3.12)}
$$

It is clear that $v_\delta(\cdot): [0, \infty)\ra L^2(\Omega)$ is
measurable. We shall prove $v_\delta(t)\in U$ for almost all $t\geq
0$. Here is the argument: When $t$ is in the set $[0,
T^*-\delta]\cap E_\delta$, we have $t+\delta\in E$. Then by (3.1),
we get $B(u^*(t+\delta), \f{\e}{2})\in U$. Since
$\|u_\delta(t)\|_{L^2(\Omega)}\leq\f{\e}{2}$ for almost all $t\geq
0$, we have
$$
\|v_\delta (t)-u^*(t+\delta)\|_{L^2(\Omega)}=\|u_\delta
(t)\|_{L^2(\Omega)}\leq \f{\e}{2}\;\;\;\mbox{for almost all}\; t
\;\mbox{in}\;  [0, T^*-\delta]\cap E_\delta,
$$
namely, $v_\delta(t)\in B(u^*(t+\delta), \f{\e}{2})$ for almost all
$t$ in the set $[0, T^*-\delta]\cap E_\delta$. Hence, $v_\delta
(t)\in U$ for almost all $t$ in the set $[0, T^*-\delta]\cap
E_\delta$. On the other hand, for almost all  $t\in [0,
T^*-\delta]\cap (E_\delta)^c$, we have $v_\delta
(t)=u^*(t+\delta)\in U$. Therefore, we have proved $v_\delta\in
\uad$.

Then, by (3.7) and (3.12), we see easily that this  control
$v_\delta$ makes the equality (3.4) hold, which leads to a
contradiction to the optimality of $T^*$ for the problem $\bf (P)$.
 Thus we have proved  $u^*(t)\in \p U$ for almost all $t\in
[0, T^*]$.

Finally, if  the control set $U$ has the  additional  property:
$\chi_{\omega} U\subset U$, then we have $\chi_\omega u^*\in \uad$.
It is clear that $y(T^*; \chi_{\omega} u^*, y_0)=y(T^*; u^*, y_0)$.
Thus, $\chi_\omega u^*$ is also  an optimal control for the problem
$\bf (P)$. Hence,  it holds that $\chi_\omega u^*(t)\in \p U$ for
almost all $t\in [0, T^*]$. This completes the proof.
\endpf

\medskip
By Theorem 3.1, we immediately get the following consequence.

\medskip
{\bf Corollary 3.2.} {\it Suppose that the control set $U$ is the
ball $ B(0,R)$ with $R>0$  and the target set $S$ is nonempty in
$L^2(\Omega)$. Let $T^*$ be the optimal time and $u^*$ be an optimal
control for the problem $\bf (P)$. Then it holds that $\|\chi_\omega
u^*(\cdot, t)\|_{L^2(\Omega)}=R$ for almost all $t\in [0, T^*]$.}

\medskip
{\bf Remark 3.3.} From the proof of Theorem 3.1, we see that  if an
admissible  control $u(\cdot, t)$ does not take its value on the
boundary of the control set  $U$   in a subset  of positive measure
in the interval $[0,T]$, where the number $T$ is such that $y(T;
u,y_0)\in S$, then there exists a "room" for us to construct another
admissible control $v$ such that the corresponding trajectory $y(t;
v,y_0)$ reaches $y(T; u, y_0)$ before the time $T$. Hence, such an
admissible control $u$ can not be optimal. This idea has been used
in [4], [11], [13] and [16]. The key point is how to use this "room"
to construct such an admissible control $v$. In this work, the null
controllability property $\bf (C)$ ( Theorem 2.1) leads us to such a
way. It was already observed in [13] that the null controllability
of the boundary controlled one-dimensional heat equation in
$(0,1)\times (0, T)$, with controls restricted on an arbitrary
subset $E\subset [0,T]$ of positive measure leads to  a bang-bang
principle of time optimal boundary controls for the one-dimensional
heat equation.
\medskip

Next, we shall use  Theorem 2.1 to derive the
 uniqueness of the optimal control for the problem
$\bf(P)$ with certain target sets and control sets.

\medskip

{\bf Theorem 3.4.} {\it Suppose that the target set $S$ is  convex
and nonempty and the control set $U$ is a closed ball. Then the
optimal control of the problem $\bf (P)$ is unique.}

\medskip

{\bf Proof.}  Let $U$ to be the closed ball $B(v_0, R)$ in
$L^2(\Omega)$, centered at  $v_0$ and of positive radius  $R$. Let
$T^*$ be the optimal time for the problem $\bf (P)$. Seeking a
contradiction, we suppose that there exist two different optimal
controls $u^*$ and $v^*$ for the problem $\bf (P)$. Then there would
exist a  subset $E_1$ of positive measure in the interval $[0,
T^*]$, such that $u^*(t)\neq v^*(t)$  for every $t\in E_1$. We first
observe that
$$
y(T^*; u^*, y_0),   y(T^*; v^*, y_0)\in S.
$$
Then we construct a control $w^*(t)$ by setting
$$
w^*(t)=\f{u^*(t)+ v^*(t)}{2} \;\; \mbox{for almost all}\; t\in [0,
\infty).
$$
It is clear that $w^*\in\uad$. Moreover, since $S$ is convex, we
have
$$
y(T^*; w^*,y_0)=\f{y(T^*; u^*,y_0)+y(T^*;v^*,y_0)}{2}\in S.
$$

On the other hand, we see that for almost all $t\in E_1$,
$$
\begin{array}{ll}
 \|w^*(t)-v_0\|^2_{L^2(\Omega)}&=2(
\|\f{u^*(t)-v_0}{2}\|^2_{L^2(\Omega)}+\|\f{v^*(t)-v_0}{2}\|^2_{L^2(\Omega)}
)-\|\f{u^*(t)-v_0}{2}-  \f{v^*(t)-v_0}{2}\|^2_{L^2(\Omega)}
\\
&=R^2-\f{1}{4}\|u^*(t)-v^*(t)\|^2_{L^2(\Omega)}\\
&<R^2.
\end{array}
$$
Thus, there exist a positive number $\e$ and a  subset $E$ of
positive measure in the set $E_1$ such that for each $t\in E$,
$\mbox{d}(w^*(t), \p B(v_0, R))\geq \e$. Then, we can use the same
argument as that in the proof of Theorem 3.1 to derive a
contradiction to
 the optimality of $T^*$. This completes the
proof.\endpf

\medskip

With regard to the existence of the time optimal controls for the
problem $\bf (P)$, we recall (See [17].) that if the target set $S$
is  closed  and convex in $L^2(\Omega)$, which contains the origin
in $L^2(\Omega)$, and if the control set $U$ is the ball $B(0,R)$
with $R>0$, then the problem $\bf (P)$ with any initial data $y_0\in
L^2(\Omega)$ has an optimal control. ( See also [14].) Thus, by
combining Corollary 3.2, Theorem 3.4 and the existence result
mentioned above, we have the following consequence.

 \medskip

 {\bf Theorem 3.5.} {\it Suppose that the target set $S$ is a closed, convex and nonempty subset,
  which  contains the origin of $L^2(\Omega)$, and the control set
 $U$ is the ball $B(0,R)$ with $R>0$.
 Then the problem $\bf (P)$ has a unique optimal control
 $u^*$ which satisfies  the bang-bang property:
 $\|\chi_\omega u^*(t)\|_{L^2(\Omega)}=R$ for almost all $t\in
 [0,T^*]$, where $T^*$ is the optimal time for the problem $\bf (P)$.}

\bigskip
{\bf Acknowledgement}. The author would like to express his
appreciation to professor Xu Zhang and Dr. K.D.Phung for their
valuable suggestions on this work.

\end{document}